\documentclass[12pt]{article}
\usepackage[english]{babel}
\usepackage[latin1]{inputenc}
\usepackage[T1]{fontenc}
\usepackage{setspace,lineno,multirow}

\usepackage{amsmath, amsfonts, amssymb, amsthm,graphicx, graphics,psfrag, caption,color}

\newtheorem{theorem}{Theorem}
\newtheorem{rk}{Remark}

\newtheorem{lem}{Lemma}

\def\calB{\mathcal{B}}
\def\N{\mathbb{N}}

\title{A Beta-splitting Model for Evolutionary Trees}
\author{Raazesh Sainudiin and Amandine V\'eber}
\date{}
\begin{document}
\maketitle

\begin{abstract}
{In this article, we construct a generalization of the Blum-Fran\c cois Beta-splitting model for evolutionary trees, which was itself inspired by Aldous' Beta-splitting model on cladograms. The novelty of our approach allows for asymmetric shares of diversification rates (or diversification `potential') between two sister species in an evolutionarily interpretable manner, as well as the addition of extinction to the model in a natural way.  We describe the incremental evolutionary construction of a tree with $n$ leaves by splitting or freezing extant lineages through the Generating, Organizing and Deleting processes.  We then give the probability of any (binary rooted) tree under this model with no extinction, at several resolutions: \emph{ranked planar trees} giving asymmetric roles to the first and second offspring species of a given species and keeping track of the order of the speciation events occurring during the creation of the tree, \emph{unranked planar trees}, \emph{ranked non-planar trees} and finally (\emph{unranked non-planar}) \emph{trees}. We also describe a continuous-time equivalent of the Generating, Organizing and Deleting processes where tree topology and branch-lengths are jointly modeled and provide code in SageMath/Python for these algorithms.}

\end{abstract}

\section{Introduction}
In the last couple of decades, many models of random evolutionary trees have been introduced and studied, as reviewed in Mooers and Heard (1997) and Morlon (2014)\nocite{mooers1997,morlon2014}. Most of them are formulated in terms of (constant or variable) individual species diversification rates mirroring the influence of particular features such as species age, trait, available niche space, etc.
In this way, they propose an evolutionary explanation for the shapes and branch lengths observed in some reconstructed real trees.
Many of these models cannot jointly model the branch lengths and the tree topologies or shapes, are quite complex to analyze and have limited identifiability (Morlon, 2014)\nocite{morlon2014}.
Note that, although we adopt here the terminology of evolutionary biology, the same kind of questions appear in other domains such as developmental biology (with cell lineage diagrams, cf. Mooers and Heard, 1997\nocite{mooers1997}, p.48) or epidemiology (Colijn and Gardy, 2014)\nocite{colijn2014}. The model developed in this paper may thus be of interest in these other contexts.

Even though the branch lengths of a phylogenetic tree give potentially precise indications on the individual diversification rates, their estimations may be subject to appreciable errors due to the difficulty of their reconstruction. On the other hand, the tree topology has a discrete nature that is somewhat easier to handle (for computation or comparison purposes, for example), and it already brings a lot of information on the phenomena shaping the clade diversity (Mooers and Heard, 1997)\nocite{mooers1997}. In particular, many works focus on the balance of a tree, measured by a diverse class of indices (for eg.~Colless index, cf. Colless, 1982,\nocite{colless1982} and Sackin index, cf. Sackin, 1972\nocite{sackin1972}).
Of course many diversification mechanisms can lead to the same phylogenetic tree balance (Jones, 2011)\nocite{jones2011} and so such indices cannot be used on their own to characterize the way the reconstructed tree was generated. However, they may be used to rule out some scenarii. For example, several papers (Mooers and Heard, 1997; Aldous, 2001; Blum and Fran\c cois, 2006)\nocite{mooers1997,aldous2001,blum2006} point at the fact that the reconstructed trees of the TreeBase database are on average much more unbalanced than expected under the most well-known model of speciation, the \emph{Yule model} (Yule, 1924)\nocite{yule1924}. In this model, every species branches into two species at the same rate (which may vary in time but remains identical for all species) and there is no extinction.
The Yule model is the best known example of an {\em evolutionarily interpretable model} of speciation due to the following three features:
\begin{itemize}
\item it is based on an incremental evolutionary construction whereby the tree grows by splitting one of the current leaf nodes which represent the set of extant lineages,
\item it can be defined jointly on the product space of tree topologies and branch lengths,
\item the distribution it induces on coarser resolutions of the tree space can be obtained.
\end{itemize}
Several models introduced in the literature are not evolutionarily interpretable in the above sense.
The main objective of this paper is to formulate an evolutionarily interpretable parametric family of models that includes the Yule model as well as many others in the literature that originally lacked evolutionary interpretability.

In Aldous (1996)\nocite{aldous1996}, Aldous introduces a one-parameter family of random cladograms, called the \emph{Beta-splitting model}. Here a cladogram is defined as a binary tree shape with a specified number of tips (or leaves) in which there is no `left' and `right' ordering of the child nodes of an internal node (in other words, the tree is \emph{non-planar} and {\em unranked} as defined below).
The leaves are labelled by the sampled species, or by $\{1,\ldots,n\}$ for simplicity.
The parameter $\beta>-2$ modulates the shape and balance of the tree produced by this model by determining the split distribution of a node subtending $m$ leaves. More precisely, Aldous' recursive construction involves a fixed $n$, the number of leaf nodes representing the extant species in a tree with at least two leaves and $\{q^\beta_n(i): 1,2,\ldots,n-1\}$, a symmetric probability distribution (i.e., $q^\beta_n(i)=q^\beta_n(n-i)$) which specifies the numbers $i$ and $n-i$ of descendants along the two branches emanating from the root node of the tree.
Once this split $(i,n-i)$ is fixed, the construction carries on recursively in the two subtrees pending from the root, with respective numbers of leaf nodes $i$ and $n-i$, and stops when all subtrees considered have only one leaf.
In the Beta-splitting model with $\beta>-2$, the split distribution $q^\beta_n$ takes the form

\begin{equation}\label{aldous}
q^\beta_n(i) = \frac{1}{a_n}\,\binom{n}{i}\int_0^1 x^{i+\beta}(1-x)^{n-i+\beta}dx
\end{equation}

for $1\leq i\leq n-1$, where $a_n$ is a normalizing factor given by

$$
a_n= \int_0^1 \big(1-x^n-(1-x)^n\big) x^\beta (1-x)^\beta dx.
$$
This \emph{Markov branching model} has now become a reference in the literature (Blum and Fran\c cois, 2006; Phillimore and Price, 2008; Jones, 2011)\nocite{jones2011,blum2006,phillimore2008}, in particular because it provides a family of random tree topologies indexed by a single parameter, which contains the most commonly used Yule tree ($\beta=0$) and Proportional to Distinguishable Arrangements (or PDA) model in which every cladogram is equally likely ($\beta=-3/2$).
The parameter $\beta$ tunes the balance of the tree, since `$\beta=-2$' corresponds to the totally unbalanced tree or comb, whereas the generated trees become more and more balanced as $\beta$ tends to infinity.
Aldous (2001)\nocite{aldous2001} also proposes a measure of the balance of a tree which has the advantage of being independent of the tree size, at least for large $n$'s: the median of the split distribution $q^\beta_n$.
This measure is used to perform maximum likelihood estimation of $\beta$ or to compare the global balance of several trees (Aldous, 2001; Blum and Fran\c cois, 2006)\nocite{aldous2001,blum2006}.

Unfortunately, Aldous, being unable to find an appropriate underlying
process (cf. Aldous, 1996, Section~4.3)\nocite{aldous1996}, in his own words, ``resort(s) to pulling a model out of thin air'' (Aldous, 1996, Section~3)\nocite{aldous1996}.
Since the number of leaf nodes has to be known before recursive splitting begins, Aldous' Beta-splitting model is not based on an incremental evolutionary construction or defined jointly on the product space of tree topologies and branch-lengths for every value of $\beta>-2$, and thus lacks evolutionary interpretability in our sense.

Subsequently, several other families of random tree topologies have been introduced, in particular Ford's alpha-model (Ford, 2005)\nocite{ford2005} in which branches are added one after another to the tree until it has the desired number of leaves. The parameter $\alpha\in [0,1]$ there serves to give a weight to each existing edge in the tree and then choose which one will be split to insert the next edge.
Ford's alpha-model also lacks evolutionary interpretability since new species can arise not just from the currently extant leaf lineages but from any ancestral lineage that is currently extinct.
Blum and Fran\c cois (2006)\nocite{blum2006} introduces an evolutionary Beta-splitting model based on ideas of Kirkpatrick and Slatkin (1993)\nocite{kirkpatrick1993}, and Aldous (1996)\nocite{aldous1996}. The idea is that the `speciation potential' is shared between the two offspring species in a random way, as may occur e.g. in the cases where speciation is influenced by available niche or geographical space that is shared between the two new species. In this model, a (rooted binary non-planar) tree is constructed incrementally by starting from a single node (the root) with speciation rate (or `potential') $1$. When this first species branches, a parameter $p_1$ is sampled in $[0,1]$ according to a Beta$(\beta+1,\beta+1)$ distribution (the definition of the Beta distribution is recalled below). Then the first offspring species is given the speciation rate $p_1$, and the second the speciation rate $1-p_1$. The next species to split is thus the first one with probability $p_1$, or the second one with probability $1-p_1$. Carrying on the construction, upon the split of a species with speciation rate $\lambda$, a new parameter $p_i$ is sampled independently of the previous ones according to the same Beta$(\beta+1,\beta+1)$ distribution, and the two sister species receive the speciation rates $\lambda p_i$ and $\lambda(1-p_i)$. Then, each species is the next one to branch with a probability equal to its speciation rate/potential.

Though the Blum-Fran\c cois and the Aldous Beta-splitting models coincide for $\beta=0$, in general they do not yield the same distribution on cladograms. See the Supplementary Material of Blum and Fran\c cois (2006)\nocite{blum2006} for a discussion of the relations between the two families of processes. Nevertheless, the principles behind the two constructions are similar and the Blum-Fran\c cois model offers an approximate evolutionary construction of Aldous' Beta-splitting model, with a slightly restricted range of parameters ($\beta>-1$ instead of $\beta>-2$). Below, we argue that the range of topologies covered by the Blum-Fran\c cois model is quite wide as well, since `$\beta=-1$' corresponds to the totally unbalanced trees while `$\beta = \infty$' corresponds to highly balanced trees. For the reasons expounded in this paragraph, we feel that this model has not yet received the attention it deserves in the phylogenetics community (or other communities as explained earlier), in particular because it is only sketchily described in Blum and Fran\c cois (2006)\nocite{blum2006}.

In this article, we extend the Blum-Fran\c cois model by allowing asymmetric Beta-distributions for the split distribution.
That is, the fraction of `speciation potential' allocated to the first offspring species is now distributed according to a Beta$(\alpha+1,\beta+1)$ distribution, for some $\alpha>-1$ and $\beta>-1$.
Of course this lack of symmetry makes sense only if we distinguish a first and second (or later `left' and `right') offspring species.
This distinction appears naturally when we think of speciation as being the creation of a new species and the continuation of the mother species, in which case the two species will not play a symmetric role and have \emph{a priori} no reason to speciate at the same rate (see e.g. Hagen et al., 2015\nocite{hagen2015}).
In the context of transmission trees in epidemiology (Sainudiin and Welch, 2015)\nocite{UCDMS20154}, the left branch keeps track of the infector and its `infection potential' while the right branch keeps track of the infectee and its `infection potential' for each infection event recorded by the internal branch.
The same distinction is true of cell lineage diagrams where the left branch can track the sister cell upon division using some measurable feature such as having more DNA damage than the sister cell along the right branch (Stewart et al., 2005)\nocite{stewart2005}.

In order to formalize more precisely how these Beta-splits create a given topology of interest, we consider four types of (rooted binary) trees:
\begin{itemize}
\item {\bf Ranked planar trees:} In this case, we distinguish the left and right child nodes of an internal node, and every internal node is labelled by an integer keeping track of the ordering in which the splits occur during the construction of the tree. Since a binary tree with $n$ leaves has $n-1$ internal nodes, the labels thus run from $1$ (the root) to $n-1$ (the last split).
\item {\bf Unranked planar trees:} Left and right child nodes are distinguished, but the internal nodes are not labelled (so that the order of the splits is not recorded).
\item {\bf Ranked non-planar trees:} In this case, the internal nodes are ranked and labelled according to the splitting order, but left and right child nodes play equivalent roles.
\item {\bf Trees:} Unranked and non-planar trees. Aldous' cladograms are such trees whose leaves are further labelled by the $n$ taxa.
\end{itemize}

Indeed, as explained above, planarity can be interesting when the two sister species do not necessarily evolve according to the same mechanisms.
The ranking of the internal nodes is a way to include some information on relative speciation times without keeping track of the full set of speciation times, see e.g. Ford et al. (2009)\nocite{ford2009}. Furthermore, various tree shape statistics are functions of the unranked non-planar trees or cladograms without leaf labels.
Explicit expressions for the probability of {\em any} tree at each of these four resolutions is not available in the literature for the Beta-splitting models of Aldous or Blum-Fran\c cois.
Thus, another contribution of this paper is the set of explicit expressions for the probability of any tree at the resolutions of ranked planar and unranked planar trees for any $\alpha$ and $\beta$ and for the probability of any tree at the resolutions of ranked non-planar and unranked non-planar trees for any $\alpha=\beta$.

We first focus on the finest of these four tree resolutions, that of ranked planar trees. We introduce the generalization of the Blum-Fran\c cois Beta-splitting model by decomposing the construction of a random ranked planar tree with $n$ leaves into two steps. First, we sample a \emph{generating sequence} $(G_i)_{i\geq 1}$, which is a realization of a sequence of independent and identically distributed random variables which, at each step $i$, will determine the choice of the next leaf to be split and the fraction of `speciation potential' allocated to the left child of that former leaf. Once this generating sequence is fixed, we define a (non random) \emph{organizing process} that turns the generating sequence into a ranked planar binary tree with the desired number of leaves. Each of these leaves is labelled by a subinterval of $[0,1]$ whose length is the speciation potential of the corresponding species.
The intervals take part in the choice of the next leaf to be split.
This construction enables us to add species extinction by a similar mechanism; thanks to a \emph{deleting process} that encodes the freezing of some leaf nodes with a given probability $\delta\in [0,1)$.
A frozen leaf can no longer evolve, and thus represents a species which is either extinct or no longer able to diversify.
See the next section for a precise description of the Generating, Organizing and Deleting processes.

Our next task is to describe the distribution on ranked planar trees corresponding to a given pair of parameters $\alpha,\beta>-1$, as well as the distribution on the three coarser tree resolutions induced by this construction.
We provide several examples in the case $\alpha=\beta$ of Blum and Fran\c cois (2006)\nocite{blum2006}, in particular to discuss the balance of the trees obtained as a function of $\beta$.
For completeness, we also propose a continuous-time process of leaf splitting and freezing such that the shape of the tree obtained after $N$ events (regardless of branch lengths) has the same distribution as that obtained through the generating, organizing and deleting processes after $N$ steps.
Finally, in the Supplementary Material we give SageMath/Python code to produce these trees at several resolutions as well as a demonstration of the code for the case of the Yule process with four leaves.

In the Supplementary Material, we also discuss a reversibility result describing how to choose a pair of sibling leaves (a \emph{cherry}) in an unranked planar tree with $n+1$ leaves created through the generating and organizing processes with $n$ steps, in such a way that the tree with $n$ leaves obtained by removing this cherry has the same law as a tree we would have obtained from the $GO$ processes with only $n-1$ steps. This last result is on \emph{sampling consistency} of our {\em evolutionarily interpretable} Beta-splitting model, which unlike Aldous' model (cf. Aldous, 1996, Section~6.3)\nocite{aldous1996}, does not naively satisfy equivalence in distribution between (i) constructing a tree with $n+1$ leaves and then removing one leaf at random and (ii) constructing a tree directly with $n$ leaves.
Our result shows that in the unranked case (planar or non planar), there is a natural but non-uniform way of choosing a terminal split to remove to obtain a tree with the same distribution as if it had been produced directly with the reduced number of leaves.

Let us end this section with some notation.
To match the standard definition of the Beta distribution, for any $\alpha,\beta>0$ we call $\calB(\alpha,\beta)$ the distribution on $[0,1]$ with density $B(\alpha,\beta)^{-1}x^{\alpha-1}(1-x)^{\beta-1}$, where

\begin{equation}\label{notation beta}
B(\alpha,\beta) := \int_0^1 x^{\alpha-1}(1-x)^{\beta-1}dx.
\end{equation}
If $\alpha=\beta$, this distribution is symmetric: if $X\sim \calB(\beta,\beta)$, then $1-X \sim \calB(\beta,\beta)$.

In all that follows, we shall consider the $\calB(\alpha+1,\beta+1)$ distribution (for $\alpha,\beta>-1$), with density proportional to $x^\alpha(1-x)^\beta$. This choice corresponds to the density used in the Aldous and Blum-Fran\c cois Beta-splitting models in the symmetric case $\alpha=\beta$.

\section{An Evolutionary Construction}
We fix $\alpha,\beta>-1$.
\subsection{The Generating Sequence}
Let $(B_1,B_2,\ldots)$ be a sequence of independent and identically distributed (i.i.d.) random variables, with the $\calB(\alpha+1,\beta+1)$ distribution. Let also $(U_1,U_2,\ldots)$ be a sequence of i.i.d. random variables with the uniform distribution on $[0,1]$, that is independent of $(B_1,B_2,\ldots)$.
Thus, each of these variables takes its values in $[0,1]$.
We call $(G_i = (U_i,B_i))_{i\in \N}$ the \emph{generating sequence}.
It will be the basis of an incremental construction of a ranked planar binary tree with $n$ leaves and $n-1$ internal nodes.

\begin{rk}
Here we use the $\calB(\alpha+1,\beta+1)$ distribution because it gives us a two-parameter family with a wide range of possible behaviours for the corresponding trees (as we shall see later). In general, we may take $(B_i)_{i\in \N}$ to be a sequence of independent and identically distributed variables with some common distribution $F$ on $[0,1]$. Even more generally, we may take a sequence with an arbitrary dependence of each $G_i$ on the previous values $(G_1,\ldots,G_{i-1})$.
\end{rk}

\subsection{The Organizing Map}
 Let us now describe the deterministic mapping that takes a realization of the generating sequence $(G_i)_{i\in \N}$ and turns it into a planar binary tree in which the internal nodes are labelled by an integer and the leaves are labelled by a subinterval of $[0,1]$. As we shall see below, the integer labels of the internal nodes will give the order in which these nodes have been split during the construction. The interval labels of the leaves will form a partition of the interval $[0,1]$ and will be used to decide which leaf is split and becomes an internal node in the next step.

Let $(g_i=(u_i,b_i))_{i\in \N}$ be a realization of the generating sequence. The organizing map $O(g)$ proceeds incrementally as follows, until the tree created has $n$ leaves. We start with a single root node, labelled by the interval $[0,1]$.
\begin{itemize}
\item Step $1$: Split the root into a left leaf labelled by $[0,b_1]$ and a right leaf labelled by $[b_1,1]$. Change the label of the root to the integer $1$.
\item Step $2$: If $u_2\in [0,b_1]$, split the left child node of the root into a left leaf and a right leaf respectively labelled by $[0,b_1b_2]$ and $[b_1b_2,b_1]$. If $u_2\in [b_1,1]$, then instead split the right child node of the root into left and right leaves with respective labels $[b_1,b_1+(1-b_1)b_2]$, $[b_1+(1-b_1)b_2,1]$. Label the former leaf that is split during this step by $2$.
\item Step $i$: Find the leaf whose interval label $[a,b]$ contains $u_i$. Change its label to the integer $i$ and split it into a left leaf with label $[a,a+(b-a)b_i]$ and a right leaf with label $[a+(b-a)b_i,b]$.
\item Stop at the end of Step $n-1$.
\end{itemize}

\begin{figure}[t]
\begin{center}
\input{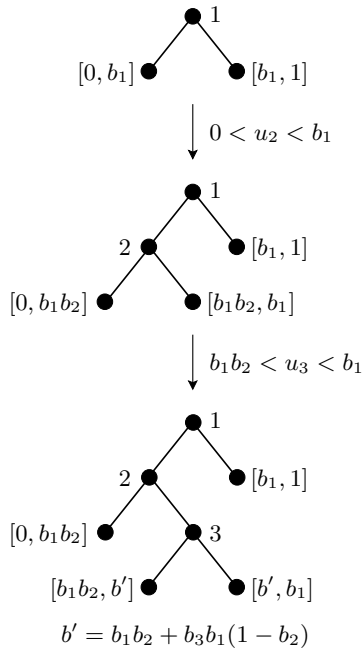}
\caption{\label{example1}An example of construction for $n=4$.}
\end{center}
\end{figure}

In words, at each step $i$ the labels of the leaves form a partition of the interval $[0,1]$. We find the next leaf to be split by checking which interval contains the corresponding $u_i$ and then $b_i$ is used to split the interval of that former leaf, say with length $\ell$, into two intervals of lengths $b_i\ell$ and $(1-b_i)\ell$. The internal node just created is then labelled by $i$ to record the order of the splits. At the end of step $i$, the tree has $i+1$ leaves, and so we stop the procedure at step $n-1$. Figure~\ref{example1} shows an example of such construction for $n=4$.

Note that once the realization of the generating sequence has been fixed, the creation of the ranked planar binary tree has no extra randomness. Below, we shall study the random tree obtained under the assumption that the generating sequence is a sequence of i.i.d. pairs $(U_i,B_i)_{i\in \N}$, where $U_i\sim \mathrm{Unif}[0,1]$ and $B_i\sim \mathcal{B}(\alpha+1,\beta+1)$.

\subsection{Generating, Organizing and Deleting Process} We can complete the organizing procedure to obtain an incremental construction of a tree with splitting (or reproduction) and freezing (or death) events. In this new process, freezing will correspond for example to a species becoming extinct (so that it cannot speciate later): such a leaf will be marked with a star and cannot be chosen to split in later steps. For this we need to augment the generating sequence to include two more coordinates, which will decide whether the next step is a split or a freezing, and which leaf is frozen in the second case.

More precisely, let $(V_1,V_2,\ldots)$ and $(D_1,D_2,\ldots)$ be two independent sequences of i.i.d.~random variables with a uniform distribution on $[0,1]$ (independent of $(G_1,G_2,\ldots)$). Let also $\delta\in [0,1)$ be a fixed number corresponding to the probability that the next event is a freezing event and not a split. We augment the generating sequence into the following sequence of quadruples $(\tilde{G}_i=(U_i,B_i,V_i,D_i))_{i\in \mathbb{N}}$.

Let now $(\tilde{g}_i=(u_i,b_i,v_i,d_i))_{i\in \N}$ be a realization of the new generating sequence. Again, we start with a single root node, labelled by the interval $[0,1]$ and proceed incrementally, until the tree created has $n$ active (i.e., not frozen) leaves or no active leaves. At each step $i$, we decide to freeze if $v_i < \delta$ and split otherwise. If $v_i < \delta$, then we freeze the leaf node whose interval label contains $d_i$ by marking it with a star (if it was already marked, then nothing changes). If $v_i \geq \delta$, we find the leaf node whose interval label contains $u_i$ as before. If the corresponding leaf is still active, we split it according to the procedure described in the organizing map. If that leaf is frozen, then the event is cancelled.
Alternatively, we can have a construction where the distribution is conditional over the currently active leaf intervals.

\begin{figure}[t]
\begin{center}
\includegraphics{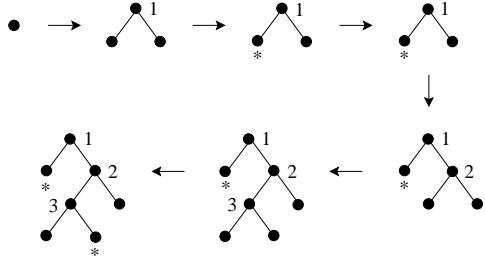}
\caption{\label{figureGOD}Example of a realization of the generating, organizing and deleting process. Here we only record the labels of the internal nodes (the split ranking) and the stars indicating a frozen leaf, but each leaf is also labelled by an interval as in the organizing process. We start with a single node. During the first step, $v_1\geq \delta$ and so the node is split and becomes labelled by $1$. Next, $v_2<\delta$ and $d_2$ belongs to the interval labelling the left leaf, so that this leaf becomes frozen. During the third step, whatever the value of $v_3$, the affected leaf chosen according to where $u_i$ or $d_i$ sits lies among the frozen leaves and so nothing happens. The next two steps are such that $v_i\geq \delta$ and the leaves chosen to split are both active. In the final step, $v_6<\delta$ and $d_6$ belongs to the interval labelling the right child leaf of node $3$, which therefore becomes frozen.}
\end{center}
\end{figure}

Figure~\ref{figureGOD} gives an example of realization of the generating, organizing and deleting process. Of course this procedure is particular in the sense that we may have chosen more general distributions for the variables $D_i$ dictating the choice of the leaf becoming frozen.

\section{Properties of the Beta-splitting evolutionary trees}
Keeping track of the generating sequence is useful to carry on the incremental construction and add new leaves to the tree. However, in most applications the object of interest is the (unlabelled) ranked planar binary tree obtained by keeping the labels of the internal nodes (giving the ranking of the splits) and by erasing the leaf nodes' interval labels whose widths give their speciation potentials that are yet to be observed.
Thus, this is the random tree of interest in this section.

\subsection{Probability of a given tree}
All trees here are rooted and binary.
First, let us give the probability of obtaining a given tree through the random generating and the non-random organizing processes.

For a given (unlabelled) ranked planar tree, and an internal node labelled by $i$, let us write $n_i^L$ (resp., $n_i^R$) for the number of internal nodes in the left (resp., right) subtree below node $i$. In particular, if node $i$ subtends two leaves, then $n_i^L=0=n_i^R$.

\smallskip
\begin{theorem}\label{th:proba}
For any unlabelled ranked planar binary tree $\tau$ with $n$ leaves, we have
\begin{align}
\mathbb{P}(\tau) &= \prod_{i=1}^{n-1} \left\{\frac{1}{B(\alpha+1,\beta+1)}\int_0^1 b_i^{n_i^L+\alpha}(1-b_i)^{n_i^R+\beta}db_i\right\} \nonumber\\
& = \prod_{i=1}^{n-1} \frac{B(n_i^L+\alpha+1,n_i^R+\beta+1)}{B(\alpha+1,\beta+1)}, \label{proba}
\end{align}
where $B(\alpha,\beta)$ was defined in \eqref{notation beta}.
\end{theorem}

\smallskip
\noindent\textbf{Proof outline.} Remember that if a leaf is labelled by an interval $[a,b]$, the probability that it is split during the $i$th step is $b-a$, the probability that the uniform random variable $U_i$ falls within $[a,b]\subset[0,1]$. If it is chosen to split, it is given label $i$ and the left and right leaves created are labelled by intervals of respective lengths $B_i(b-a)$ and $(1-B_i)(b-a)$. Then these intervals may split later, but into intervals of lengths that are always proportional to $B_i$ or $1-B_i$ (respectively). Now the probability of the tree $\tau$ is the product of the $n-1$ probabilities of choosing a given leaf to split at each step, each of which is equal to the length of the interval labeling that leaf. As a consequence, each split occurring in the left subtree below node $i$ brings in another $B_i$ in the product, or another $1-B_i$ if the split occurs in the right subtree below node $i$. Averaging over the possible values of the $B_i$'s, which are independent $\mathcal{B}(\alpha+1,\beta+1)$ random variables, yields the result. \hfill $\Box$

\begin{rk}
This construction is different from Aldous' interpretation in terms of splitting intervals that starts by uniformly scattering the given $n$ leaf nodes as `particles' on the unit interval and splitting the interval at a random point with density $f$.
This splitting is repeated recursively on sub-intervals exactly as we do, i.e.~splitting each interval $[a,b]$ at a point $a+X(b-a)$ where the $X$'s are independent with density $f$.
Splitting stops when each subinterval contains only one leaf particle while splits that result in one of the intervals being empty (without any leaf particles in it) are not allowed.
See the Supplementary Material of Blum and Fran\c cois (2006)\nocite{blum2006} for a discussion on the relation between Aldous' Beta-splitting model and this incremental construction.
\end{rk}

\subsection{Examples}
In all the examples given below, we focus on the symmetric case $\alpha=\beta$. Some of the formulae given below are easily generalized to the case $\alpha\neq \beta$.

The most important example is the case $\beta=0$, which corresponds to the Yule model of pure births that is used in many models of phylogenies.

Recall that $B(\alpha,\beta)$ is related to the Gamma function $\Gamma$ by the equality

\begin{equation}
B(\alpha,\beta)=\frac{\Gamma(\alpha)\Gamma(\beta)}{\Gamma(\alpha+\beta)},\qquad \alpha,\beta>0,
\end{equation}
and that $\Gamma(\beta)=(\beta-1)!=(\beta-1)(\beta-2)\cdots 2\cdot 1$ if $\beta\in \N$. Using \eqref{proba} with $\alpha=\beta = 0$, we have

\begin{equation}\label{beta0}
\mathbb{P}(\tau)= \prod_{i=1}^{n-1} \frac{n_i^L!n_i^R!}{(n_i^L+n_i^R+1)!} = \frac{1}{(n-1)!},
\end{equation}
where the second equality is obtained by observing that $n_i^L+n_i^R+1$ is the number of internal nodes of the tree rooted at node $i$, which is the left or the right subtree below the mother node of $i$. Hence, each term $n_i^L!$ in the numerator of the product cancels with the term in the denominator that corresponds to the left child node of $i$, except if $n_i^L=0$ and the left child node of $i$ is a leaf. But in this case, $0!=1$ by convention. The same holds true for each of the $n_i^R!$. Likewise, the terms in the denominator which are not compensated by some term in the numerator are those corresponding to internal nodes having no mother nodes. But the only such node is the root ($i=1$), with $n_1^L+n_1^R+1=n-1$. This gives us the result.

\begin{rk}
This construction is very different from the standard evolutionary construction of the Yule tree, in which the next leaf to split is chosen uniformly at random among the current set of leaves.
Here the choice of the next split is dictated by the lengths of the intervals labeling the current leaves, which will all be distinct with probability one.
However, averaging over the law of the generating sequence (when $\alpha=\beta=0$) yields the same distribution on ranked planar binary trees.
\end{rk}

Using the above property of the Gamma function, we can also give explicit values for the probability of a tree when $\alpha=\beta$ is a non-negative integer: if $\beta\in\N\cup \{0\}$, then

\begin{equation}\label{beta integer}
\mathbb{P}(\tau) = \prod_{i=1}^{n-1}\frac{(n_i^L+\beta)!(n_i^R+\beta)!(2\beta+1)!}{(n_i^L+n_i^R+2\beta+1)!(\beta!)^2}.
\end{equation}

Thirdly, when $\beta$ is a nonnegative integer, we have

$$
\Gamma(\beta+1/2)= \frac{(2\beta)!}{2^{2\beta}\beta!}\, \sqrt{\pi}.
$$
As a consequence, another example in which the probability of a tree has an explicit form is the case where $\alpha=\beta=b-1/2$, with $b\in \N\cup \{0\}$:

$$
\mathbb{P}(\tau) = \prod_{i=1}^{n-1} \frac{(2n_i^L+2b)!(2n_i^R+2b)!(b!)^2}{4^{n_i^L+n_i^R}(n_i^L+b)!(n_i^R+b)!(n_i^L+n_i^R+2b)!(2b)!}.
$$
To our knowledge, the cases $\alpha=\beta\in \N$ and $\beta +1/2\in \N\cup \{0\}$ correspond to no well-studied models of trees.

\begin{figure}[t]
\begin{center}
\input{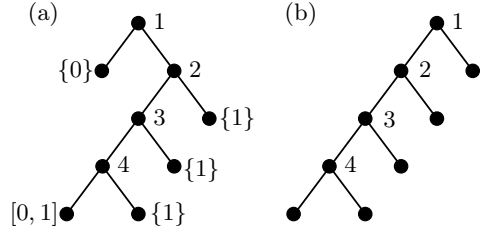}
\caption{\label{example3}(a) An example of realization of a tree corresponding to the limiting case $\beta=-1$, and (b) the comb which is the only possible non-planar tree that can be generated in this case.}
\end{center}
\end{figure}

To see how the global shape of the tree (and in particular its balance) evolves as $\beta$ goes from $-1$ to $+\infty$, let us consider the two extreme cases. The corresponding processes cannot be defined directly as $-1$ and $+\infty$ lie out of the range of the possible $\beta$'s, but we can capture the essence of the resulting (random) tree by taking limits in $\beta$. First, as $\beta\rightarrow -1$, the $\calB(\beta+1,\beta+1)$ distribution gives more and more weight to the boundaries $0$ and $1$. In the limit, the random variables $B_i$ should then take the values $0$ or $1$, each with probability $1/2$. In this case, the root is first split into a leaf with label $[0,1]$ and another leaf with label $\{0\}$ or $\{1\}$ (i.e., an interval reduced to a single point). The leaf that receives the label $[0,1]$ is the left one with probability $1/2$. Next, the uniform random variable $U_2$ belongs to the interval $[0,1]$ with probability one, so that the leaf labelled by $[0,1]$ is necessarily that chosen to split. Again, it is split into two leaves with labels $[0,1]$ and $\{0\}$ or $\{1\}$, implying that the next leaf to split is that inheriting the full interval $[0,1]$ with probability one. The reasoning can be carried on until step $n-1$. Hence, morally the tree corresponding to $\alpha=\beta=-1$ is a fully unbalanced tree, with a single backbone from which the $n$ leaves are hanging. The backbone is extended at each step by choosing one of the two leaves created in the previous step, each with probability $1/2$. See Figure~\ref{example3} for an example with $n=5$.

Let us now consider the limit $\beta\rightarrow +\infty$. Using \eqref{beta integer} and the fact that $(b+i)!/b! \sim b^i$ as $b\rightarrow \infty$ (meaning that the ratio of both terms tends to $1$), we can pass to the limit $\beta\rightarrow \infty$ and obtain that

$$
\lim_{\beta\rightarrow +\infty}\mathbb{P}(\tau) = \prod_{i=1}^{n-1}\frac{1}{2^{n_i^L+n_i^R}}.
$$
Because in the balanced trees the internal nodes below a given node are equally split between the left and right subtrees hanging from that node, the sum $n_i^L+n_i^R$ decreases with $i$ faster than in more unbalanced trees. This means that for very large $\beta$'s, approximately balanced tree will have much higher probabilities than unbalanced ones. For instance, any of the fully unbalanced trees $\tau_u$ will have probability

$$
\mathbb{P}(\tau_u)=1/[2^{n-2}\cdot 2^{n-3}\cdots 2\cdot 1] = 2^{-(n-1)(n-2)/2}.
$$
On the other hand, if $n = 2^N$ is a power of 2, the probability of any of the fully balanced trees $\tau_b$ is equal to

\begin{align*}
\mathbb{P}(\tau_b)&= 1/\bigg[ 2^{n-2}.(2^{\frac{n}{2}-2})^2\cdot(2^{\frac{n}{4}-2})^4 \cdots (2^2)^{2^{N-2}} \cdot (2^0)^{2^{N-1}}\bigg]\\
& = \prod_{k=0}^{N-1}\big(2^{\frac{n}{2^k}-2}\big)^{-2^k} = \prod_{k=0}^{N-1}2^{2^{k+1}-n} = 2^{-n(N-2)-2}.
\end{align*}
Indeed, any subtree pending from a node at level $k \in \{0,\ldots,N-1\}$ (level $0$ being that of the root, level $N=\log_2(n)$ that of all the leaves) has $n/2^k$ leaves, and so $(n/2^k)-2$ internal nodes below its root. Furthermore, there are $2^k$ internal nodes at level $k$. Together with \eqref{proba}, this gives us the result.

Using the results obtained in the next section, we can further compute the probability of producing a fully unbalanced (unranked nonplanar) tree as being equal to
$$
\mathbb{P}(t_u) = 2^{n-2}\,\mathbb{P}(\tau_u) = \left\{\begin{array}{ll}
1 & \quad \hbox{if }\beta=-1,\vspace{0.1cm}\\
\frac{2^{n-2}}{(n-1)!} & \quad \hbox{if }\beta=0,\vspace{0.1cm}\\
2^{-(n-2)(n-3)/2} & \quad \hbox{if }\beta=+\infty.
\end{array}\right.
$$
Likewise, the probability of producing a fully balanced tree with $n=2^N$ tips is given by
$$
\mathbb{P}(t_b) = \frac{(n-1)!}{\prod_{k=0}^{N-1}\left(\frac{n}{2^k}-1\right)^{2^k}}\, \mathbb{P}(\tau_b)
$$
Table~\ref{probabilities} gives a few examples of these probabilities for different values of $n=2^N$ and $\beta = -1,0,+\infty$.

\begin{table}
\centering
\begin{tabular}{|c|c|c|c|c|c|}
  \hline
    &  $n$& $4$ & 8 & 32 & 1024 \\ \hline
  \multirow{2}{*}{$\beta=-1$} & Comb & 1&1&1&1 \\
    & Balanced & 0&0&0&0 \\\hline
  \multirow{2}{*}{$\beta=0$ } & Comb & $0.667$ & $1.27e^{-2}$ & $1.31e^{-25}$ & $8.49e^{-2330}$ \\
    & Balanced & $0.333$ & $1.59e^{-2}$ & $9.10e^{-12}$ & $1.04e^{-417}$ \\ \hline
  \multirow{2}{*}{$\beta=+\infty$} & Comb & $0.5$ & $3.05e^{-5}$ & $1.13e^{-131}$ & $2.09e^{-157057}$\\
    & Balanced & $0.5$ & $7.81e^{-2}$ & $2.36e^{-7}$ & $1.26e^{-247}$ \\
 \hline
\end{tabular}
\caption{\label{probabilities} Probability of sampling a comb tree or a fully balanced tree for different values of $n$ and $\beta$. As explained in the text, larger values of $\beta$ correspond to higher probabilities of sampling a balanced tree.}
\end{table}

Finally, we have shown that the family of Beta-splitting trees defined in Blum and Fran\c cois (2006)\nocite{blum2006} and generalized in this article includes a one-parameter family containing the classical Yule (ranked planar) tree. For small $\beta$'s (close to $-1$), the corresponding trees are unbalanced with high probability, whereas for large $\beta$'s the tree distribution is concentrated on balanced trees. The family of Beta-splitting trees indexed by $\alpha,\beta>-1$ thus covers a very wide range of possible topologies.

\section{Other tree resolutions}
Recall that a tree in this paper is always rooted and binary.
Up to now we have focused on ranked planar trees with $n$ leaves that keep records of the order in which splits occur and give an asymmetrical role to the left and right child nodes of an internal node.
These $(n-1)!$ many trees are in bijective correspondence with permutations of $\{1,\ldots,n-1\}$ through the {\em increasing binary tree-lifting} operation (see Flajolet and Sedgewick, 2009, Ex.~17, p.~132)\nocite{Flajolet}.
However, we may be interested in coarser resolutions of the trees generated by our Beta-splitting procedure, especially those resolutions of interest to evolutionary biologists.
To the best of our knowledge, explicit formulae for the probability of any tree at these resolutions are not available in the literature as a function of $\alpha$ and $\beta$ (even for the symmetric case when $\alpha=\beta$).
This is because the cardinality of the inverse image from a fine to a coarser tree resolution needs to be computed for any tree in the coarser resolution.
Such probabilities can be directly useful in simulation-intensive inference.

\subsection{Probability of unranked planar trees}
Here we keep the lack of symmetry between the child nodes, but do not record the order of the splits. As explained in the introduction, this may be of interest for example if we assume that there is a lack of symmetry between the two species created during a speciation event, say due to one species being the ancestor and other being the descendant, but we do not want to reconstruct the temporal order in which the speciation events occurred.
In the context of transmission trees, we may only be interested in the infector-infectee relation for each transmission event and not in the ranking of transmission events given by their relative temporal order.

Since we do not label the internal nodes, let $\mathcal{I}(\mathbf{t})$ denote the set of all internal nodes of a planar tree $\mathbf{t}$ and let us extend the notation $n_i^L$ and $n_i^R$, $i\in \mathcal{I}$, for the number of internal nodes in the left and right subtrees below node $i$ to this unlabelled case. The probability of obtaining a given (unranked) planar binary tree $\mathbf{t}$ through the Beta-splitting generating and organizing processes is given by the following lemma.

\begin{lem}\label{proba unranked}
Let $\mathbf{t}$ be a planar binary tree. We have
\begin{align*}
\mathbb{P}(\mathbf{t}) & = \prod_{i \in \mathcal{I}(\mathbf{t})}\binom{n_i^L+n_i^R}{n_i^L}\prod_{i \in \mathcal{I}(\mathbf{t})}\frac{B(n_i^L+\alpha+1,n_i^R+\beta+1)}{B(\alpha+1,\beta+1)}\\
& = (n-1)!\prod_{i \in \mathcal{I}(\mathbf{t})}\frac{B(n_i^L+\alpha+1,n_i^R+\beta+1)}{(n_i^L+n_i^R+1)B(\alpha+1,\beta+1)}.
\end{align*}
\end{lem}
Indeed, recall that the second product in the right-hand side of the first equality above is the probability of a given ranked planar tree corresponding to the unranked tree $\mathbf{t}$. Since it does not depend on the ranking, there remains to count the number of ranked trees whose unranking gives $\mathbf{t}$.
Now, to rank the internal nodes of $\mathbf{t}$, at each split we have to decide which of the remaining integer labels go to the left or to the right subtree below the corresponding node.
This gives us $\mathrm{Binomial}(n_i^L+n_i^R,n_i^L)$ choices, hence the first product term in $\mathbb{P}(\mathbf{t})$.
This product of binomial coefficients is called the shape functional (Dobrow and Fill, 1995)\nocite{Dobrow1995}, the Catalan coefficient (Sainudiin, 2012)\nocite{CatalanCoeff2012} and is the solution to an enumerative combinatorial exercise (Stanley, 1997, Ch.~3, Ex.~1.b, p.~312)\nocite{Stanley1997}.

As in the case $\alpha=\beta=0$ (see the derivation of \eqref{beta0}), the simplification leading to the last equality comes from the fact that $n_i^L+n_i^R+1$ is the number of internal nodes of the subtree rooted at node $i$, so that most factorial terms cancel out in the product over $\mathcal{I}(\mathbf{t})$.

\subsection{Probability of ranked non-planar trees}
In this case we keep the ranking but give a symmetric role to the left and right child nodes of an internal node.
These trees are termed {\em evolutionary relationships} by Tajima (Tajima, 1983)\nocite{tajima1983} who shows that there are $2^{n-1-c(\mathfrak{t})}$ ranked planar trees for a given ranked non-planar tree $\mathfrak{t}$, where $c(\mathfrak{t})$ is the number of cherry nodes, i.e.~sub-terminal nodes with two child nodes. For a quick justification of Tajima's result, suppose we want to turn the ranked non-planar tree $\mathfrak{t}$ into a planar tree. For each of the $n-1$ internal nodes of $\mathfrak{t}$, there are 2 choices for the child node that is said to be `left' except if they are both leaves (i.e., the internal node is a cherry node). Indeed, in this case they carry no ranking that would make them distinguishable.

Since the probability of a ranked non-planar tree does not depend on the planarity, provided $\alpha=\beta$, the probability of $\mathfrak{t}$ is simply the product of the probability of a corresponding ranked planar tree, times the number of ranked planar trees corresponding to $\mathfrak{t}$. That is:

\begin{align*}
\mathbb{P}(\mathfrak{t}) & = 2^{n-1-c(\mathfrak{t})} \prod_{i=1}^{n-1}\frac{B(n_i^L+\beta+1,n_i^R+\beta+1)}{B(\beta+1,\beta+1)}.
\end{align*}
Note that when $\alpha \neq \beta$ we need to sum over all $2^{n-1-c(\mathfrak{t})}$ ranked planar trees that map to the ranked non-planar tree $\mathfrak{t}$ (since in this case $B(n_i^L+\alpha+1,n_i^R+\beta+1) \neq B(n_i^R+\alpha+1,n_i^L+\beta+1)$), and this may not be computationally feasible for large $n$.

\subsection{Probability of trees}
This is the case of (rooted binary) unranked non-planar trees or simply trees (also called phylogenetic tree shapes).
There are $2^{n-1-s(t)}$ unranked planar trees that correspond to a tree $t$, where $s(t)$ is the number of internal nodes of $t$ that have isomorphic left and right subtrees.
See Sainudiin et al. (2015)\nocite{UCDMS20152} for a proof by induction.
And all these unranked planar trees that correspond to $t$ have the same probability provided $\alpha=\beta$.
One can intuitively understand this by noting that there are two unranked planar embeddings for each internal node of $t$ that does not have isomorphic subtrees on its left and right descendant nodes.
Thus, the probability of a tree $t$ if $\alpha=\beta$ is:

\begin{align*}
\mathbb{P}(t) & = 2^{n-1-s(t)}
(n-1)!\prod_{i \in \mathcal{I}(t)}\frac{B(n_i^L+\beta+1,n_i^R+\beta+1)}{(n_i^L+n_i^R+1)B(\beta+1,\beta+1)}.
\end{align*}
Once again if $\alpha \neq \beta$ we need to sum over all $2^{n-1-s(t)}$ unranked planar trees that map to the (unranked non-planar) tree $t$ and this may not be computationally feasible for large $n$.

\section{Continuous-time Process}
Up to now we have described the generating-organizing-deleting process in discrete time over ranked planar trees and given the probabilities over various equivalence classes of trees. However, one may want to formulate a continuous-time version of this process in order to have a more precise description of the evolutionary relationships including how much time elapsed between two speciation or extinction events.

To do so, recall that we fix two parameters $\alpha,\beta>-1$ characterizing the way in which leaf intervals are split, and the probability $\delta\in [0,1)$ of a freezing during the next event (if $\delta=0$, only splits occur). Let us also fix a rate $\lambda>0$ of events. One way to formulate a continuous-time generating-organizing-deleting process is the following: suppose the current interval length of the $j$-th active leaf is $L_j$. Then each active leaf $j$ splits at rate $\lambda(1-\delta)L_j$ or becomes frozen at rate $\lambda\delta L_j$. When a split occurs, the internal node created during the event is labelled by the first integer $N$ larger than all integer labels in the current tree, and the two leaves created are labelled by intervals that are obtained by splitting the (former) interval label of node $N$ using $b_N$ (that is, if that interval is $[a,b]$, the new leaves are labelled by $[a,a+b_N(b-a)]$ and $[a+b_N(b-a),b]$ as before).

\begin{lem}~\label{lem:continuous}
The discrete tree embedded in the continuous-time planar ranked tree stopped after the $n$-th event has the same law as the (ranked planar) tree obtained from the generating, organizing and deleting process stopped after the $n$-th \emph{effective} event. By effective event, we mean an event affecting an active leaf and therefore leading to a change in the tree.
\end{lem}

\smallskip
\noindent\textbf{Proof.} Let us write $L_{\mathrm{act}}$ for the sum over all active leaves of their interval lengths $L_j$, and $I_{\mathrm{act}}$ for the union of the corresponding intervals. Hence, $1-L_{\mathrm{act}}$ is the length of the set $[0,1]\setminus I_{\mathrm{act}}$ corresponding to all frozen leaves. Let us check for both random trees that, conditionally on the current state of the process:
\begin{itemize}
\item[$(i)$] The next (effective) event is a split with probability $1-\delta$ or a freezing with probability $\delta$.
\item[$(ii)$] If it is a split, then the probability that leaf $j$ is chosen to split is $L_j/L_{\mathrm{act}}$.
\item[$(iii)$] If it is a freezing, then the probability that leaf $j$ is chosen to freeze is also $L_j/L_{\mathrm{act}}$.
\end{itemize}
The result is an easy consequence of the construction of continuous-time jump processes for the tree embedded in the continuous-time procedure (in essence, if we have a countable collection of events such that event $i$ happens at rate $\mu_i$, then the first event to occur is the $j$-th one with probability $\mu_j/\sum_i \mu_i$).

For the tree constructed from the discrete process `restricted' to the effective events, observe first that the coordinates $\{u_i,\, d_i,\ i\in \N\}$ of the organizing process are recorded in the ranked planar tree only through the choices of the next leaf to be affected. Also, the coordinates $\{v_i,\, i\in \N\}$ appear only through the types of the next events to occur. As a consequence, the law of the tree emanating from this construction depends only on the probabilities that each of these quantities belongs to a given set, conditionally on the fact that the leaf chosen is active. That is, the probability that the next effective event is a freezing is

$$
\mathbb{P}(V_i<\delta\, |\, D_i \in I_{\mathrm{act}}) = \frac{\mathbb{P}(V_i<\delta) \mathbb{P}(D_i \in I_{\mathrm{act}})}{\mathbb{P}(D_i \in I_{\mathrm{act}})} = \delta
$$
since $V_i$ and $D_i$ are independent. Likewise, the probability of the next effective event being a split is $1-\delta$, which proves $(i)$. Next, conditionally on the next event being a split, the probability that leaf $j$ is chosen is (again by the independence of $U_i$ and $V_i$)

$$
\mathbb{P}(U_i\in I_j\, |\, V_i>\delta ,\, U_i \in I_{\mathrm{act}}) = \mathbb{P}(U_i\in I_j\, |\, U_i \in I_{\mathrm{act}}) = \frac{L_j}{L_{\mathrm{act}}}.
$$
This proves $(ii)$, and $(iii)$ can be obtained in the same way. Points $(i)$, $(ii)$ and $(iii)$ then enable us to conclude that the topologies and rankings of both trees have the same law.

There remains to show that, conditionally on the topology and ranking, the interval labels of the leaves are identical in distribution. Notice that they are not \emph{a priori} identical with probability one since the continuous-time construction uses the variables $(B_1,B_2,\ldots,B_n)$ whereas the discrete-time construction uses the variables $(B_{i_1},B_{i_2},\ldots,B_{i_n})$, where $i_1<i_2< \cdots <i_n$ are the random indices of the effective events. But the event that $i_j= k$ depends only on the generating sequence $(\tilde{G}_i)_{1\leq i\leq k-1}$, since it depends only on the current length of all active leaves. Consequently, $i_j$ and $B_{i_j}$ are independent random variables (recall that all components of $(\tilde{G}_i)$ are independent of each other). A simple argument then shows that the law of $(B_{i_1},B_{i_2},\ldots,B_{i_n})$ is the same as the law of $(B_1,B_2,\ldots,B_n)$, which in turn guarantees that the interval labels of the leaves are also equal in law. Lemma~\ref{lem:continuous} is proved. \hfill $\Box$

\bigskip
\noindent {\bf Algorithm.}
The code developed for this work is given in the Supplementary Material. It is also publicly shared at \texttt{https://cloud.sagemath.com/projects/2c5f7f68-e689-4c70-a4b4-}\\
\texttt{5b5d4dc4f93f/files/2015-10-27-082849.sagews}.

\bigskip
\noindent {\bf Competing interests.}
We have no competing interests.

\bigskip
\noindent {\bf Authors' contributions.} R.S. and A.V. designed the model, studied its properties and drafted the manuscript. R.S. coded the algorithm. All authors gave final approval for publication.

\bigskip
\noindent {\bf Acknowledgments.} R.S.~thanks Mike Steel for posing the problem and R.S. and A.V.~thank H\'el\`ene Morlon for generous introductions to models of diversification. They also thank the referees for their careful reading of the manuscript. R.S.~was partly supported by a Sabbatical Grant from College of Engineering (University of Canterbury, NZ), External Consulting for Wynyard Group (Christchurch, NZ), and a Visiting Scholarship at Department of Mathematics (Cornell University, USA). R.S.~and A.V.~were supported in part by the chaire Mod\'elisation Math\'ematique et Biodiversit\'e of Veolia Environnement-\'Ecole Polytechnique-Museum National d'Histoire Naturelle-Fondation X.

\newpage
\section*{Supplementary Material}
\subsection*{Algorithm}
This code is publicly shared at \texttt{https://cloud.sagemath.com/projects/2c5f7f68-e689-4c70-}
\texttt{a4b4-5b5d4dc4f93f/files/2015-10-27-082849.sagews}.
The code was mainly used to aid intuition during this study and is not written to be efficient for large scale simulation studies.
The core Algorithms for the generating and organizing processes are presented as SageMath/python code instead of pseudo-code in order to communicate the Algorithms used in this study in a more concrete and reproducible manner.
This also allows the reader to perform computational experiments in SageMath/python immediately to further extend this work.

\medskip
{\small
The function \texttt{split01ScaledCD} takes the interval $I$ and splits it into 2 intervals of lengths $x|I|$ and $(1-x)|I|$ (with $x\in [0,1]$).

\smallskip
\begin{verbatim}
def split01ScaledCD(x,I):
    '''x \in [0,1], c=I[0] < d=I[1]'''
    c=I[0]; d=I[1];
    return [[c,c+(d-c)*x],[c+(d-c)*x,d]];
\end{verbatim}

\smallskip
\noindent The function \texttt{SplittingPermutation} translates a sequence of $n$ real numbers (our splitting points, later) into a permutation of $[n]=\{1,\ldots,n\}$ by returning a list of $n$ integers such that the $i$-th element of the list is the index of the $i$-th smallest number in the initial sequence. For example, \texttt{SplittingPermutation}($[1.1, 10, -1, 2.5]$)= $[3,1,4,2]$.

\smallskip
\begin{verbatim}
def SplittingPermutation(splitpointsequence):
    '''return the permutation of [n] given by the map from
    the sequence of n real numbers in the list
    splitpointsequence to an ordering by indices in [n]'''
    sss=sorted(splitpointsequence)
    return tuple([splitpointsequence.index(i)+1 for i in sss])
\end{verbatim}

\smallskip
In the function \texttt{MakePartitionAndTree}, we construct $m$ samples of a tree with $n$ splits (or $n+1$ leaves), using the $\calB(a+1,b+1)$ distribution of the coordinates $B_i$ of the generating sequence.
We first obtain the sequence of points in $[0,1]$ where the splits occur according to the generating sequence and store them in {\tt SplitPoints}.
Then we use the standard {\tt binary\_search\_insert} method to obtain a binary search tree that organizes the points in {\tt SplitPoints} into an unranked planar binary tree.
Recall that the split points are inserted from the root of the tree such that the new point that is smaller/larger than the point at the root node descends into the left/right subtree of the root and recursively takes left/right subtree depending on whether it is lesser or greater than the next point it encounters at an internal node that is already stored in the tree.
Each of these trees is recorded at several resolutions: the ranked planar trees (recorded in the list \texttt{SplittingPermutationSamples}, unranked planar trees (recorded in \texttt{PlanarShapeSamples}) and unranked non-planar trees (recorded in \texttt{PhyloShapeSamples}).

For the finest resolution of ranked planar trees, we use the bijection between ranked planar trees with $n+1$ leaves and permutations of $[n]$ (note that the cardinality of each set is $n!$).
This bijection is detailed in Flajolet and Sedgewick (2009), Ex.~17, p.~132,\nocite{Flajolet} and follows this idea.
Say the permutation of $[n]$ of which we want to draw the tree is $[i_1,i_2,\ldots,i_n]$.
We first construct the planar skeleton of the ranked internal nodes incrementally by starting with a single node.
To place the second node, check whether $i_2<i_1$, in which case the left child of the root becomes node $2$, or whether $i_2>i_1$ and the right child node of the root becomes node $2$.
To place the third node, find whether it goes to the left or to the right of the root by checking whether $i_3<i_1$ or $i_3>i_1$.
Once this is decided, if the second and third nodes are on the same side of the root, then compare $i_3$ to $i_2$ to decide whether node $3$ should be the left or right child of node $2$.
Proceeding in the same way for the other terms of the permutation, we construct a planar ranked skeleton with $n$ nodes.
There remains to attach the $n+1$ leaves to the terminal nodes of the skeleton to obtain a ranked planar tree with $n$ splits.
For ease of representation, the output corresponding to the ranked planar trees is thus the list \texttt{SplittingPermutationSamples} recording the permutations corresponding to the $m$ trees.

\smallskip
\begin{verbatim}
def MakePartitionAndTree(n,m,a,b):
    '''This creates m independent trees with n+1 leaves and alpha=a, beta=b
       n >= 1, where n+1 is the number of leaves
       m is the number of replicates
       a,b>-1, where (a+1,b+1) are the parameters of the beta distribution'''
    PlanarShapeSamples=[];
    PhyloShapeSamples=[];
    SplittingPermutationSamples=[];
    BetaD = RealDistribution('beta', [a+1, b+1],seed=0)
    show(BetaD.plot(xmin=0,xmax=1),figsize=[7,2]);
    print '\n';
    for reps in range(m):
        # generating i.i.d. samples from BetaD
        B=[BetaD.get_random_element() for _ in range(n)]
        #initialize
        SplitPoints=[B[0]]# keep order of split points
        Splits=split01ScaledCD(B[0],[0,1])

        #iterate
        for i in range(1,n):
            Widths=[x[1]-x[0] for x in Splits];
            W = GeneralDiscreteDistribution(Widths);
            nextSplitI=Splits[W.get_random_element()];
            Splits.remove(nextSplitI);
            NewLeaves=split01ScaledCD(B[i],[nextSplitI[0], nextSplitI[1]]);
            # find the split point between the new leaves
            RescaledG=NewLeaves[0][1];
            SplitPoints.append(RescaledG);
            Splits.append(NewLeaves[0]);
            Splits.append(NewLeaves[1]);
        SplittingPermutationSamples.append(SplittingPermutation(SplitPoints));
        # insert the split points into the tree
        t = LabelledBinaryTree(None)
        for i in range(0,n):
            t = t.binary_search_insert(SplitPoints[i]);
        sh=t.shape();
        PlanarShapeSamples.append(sh);
        PhyloShapeSamples.append(Graph(sh.to_undirected_graph(with_leaves=True),
                                                                immutable=True));
    return (PlanarShapeSamples,PhyloShapeSamples,SplittingPermutationSamples);
\end{verbatim}

\smallskip
\noindent We also provide some additional functions computing the probabilities of a given tree at a particular resolution under the Beta-splitting model.

\smallskip
\begin{verbatim}
def splitsSequence(T):
    '''return a list of tuples (left,right) split sizes at each split node'''
    l = []
    T.post_order_traversal(lambda node:
       l.append((node[0].node_number(),node[1].node_number())))
    return l

def isIso(N):
    '''does node N of binary tree have the same left and right subtree shapes
    (are left and right subtrees of node N in tree isomorphic)'''
    L=Graph(N[0].canonical_labelling().shape().to_undirected_graph(with_leaves=True),
                                          immutable=True)
    R=Graph(N[1].canonical_labelling().shape().to_undirected_graph(with_leaves=True),
                                          immutable=True)
    return 1 if L==R else 0

def numIso(T):
    '''number of internal nodes that have isomorphic left and right sub-trees'''
    l = []
    T.post_order_traversal(lambda node:l.append(isIso(node)))
    return sum(l)

def prob_RPT(T,a,b):
    '''probability of ranked planar tree T under beta-splitting model
       a,b>-1, where (a+1,b+1) are the parameters of the beta distribution'''
    # non-cherry splits
    ncspS=filter(lambda x: x!=(0,0), splitsSequence(T))
    return prod(map(lambda x:beta(x[0]+a+1,x[1]+b+1)/beta(a+1,b+1),ncspS))

def prob_PT(T,a,b):
    '''probability of planar tree T under beta-splitting model
       a,b>-1, where (a+1,b+1) are the parameters of the beta distribution'''
    # non-cherry splits
    ncspS=filter(lambda x: x!=(0,0), splitsSequence(T))
    return prod(map(lambda x: binomial(x[0]+x[1],x[1])*
            beta(x[0]+a+1,x[1]+b+1)/beta(a+1,b+1),ncspS))

def prob_RT(T,a,b):
    '''probability of ranked (nonplar) tree T under beta-splitting model
       a,b>-1, where (a+1,b+1) are the parameters of the beta distribution'''
    assert(a==b)
    spS=splitsSequence(T)
    numSplits=len(spS)
    # non-cherry splits
    ncspS=filter(lambda x: x!=(0,0), spS)
    numCherries=numSplits-len(ncspS)
    probRPT = prod(map(lambda x:beta(x[0]+a+1,x[1]+b+1)/beta(a+1,b+1), ncspS))
    return 2^(numSplits-numCherries)*probRPT

def prob_T(T,a,b):
    '''probability of tree T (phylo tree shape) under beta-splitting model
       a,b>-1, where (a+1,b+1) are the parameters of the beta distribution'''
    assert(a==b)
    spS=splitsSequence(T)
    numSplits=len(spS)
    # non-cherry splits
    ncspS=filter(lambda x: x!=(0,0), spS)
    probPT = prod(map(lambda x: binomial(x[0]+x[1],x[1])*
             beta(x[0]+a+1,x[1]+b+1)/beta(a+1,b+1),ncspS))
    numIsoSplits=numIso(T)
    return 2^(numSplits-numIsoSplits)*probPT

def stats_probs_Tree(T,a,b):
    '''probability of various resolutions of tree T under beta-splitting model
       a,b>-1, where (a+1,b+1) are the parameters of the beta distribution'''
    spS=splitsSequence(T)
    numSplits=len(spS)
    # non-cherry splits
    ncspS=filter(lambda x: x!=(0,0), spS)
    numCherries=numSplits-len(ncspS)
    probRPT = prod(map(lambda x:beta(x[0]+a+1,x[1]+b+1)/beta(a+1,b+1), ncspS))
    catCoeff = prod(map(lambda x:binomial(x[0]+x[1],x[1]),ncspS))
    # prob of (non-ranked) planar tree
    probPT = catCoeff * probRPT
    probRT = 2^(numSplits-numCherries)*probRPT
    numIsoSplits=numIso(T)
    probT=2^(numSplits-numIsoSplits)*probPT
    return (numSplits,numIsoSplits,numCherries,catCoeff,probRPT,probPT,probRT,probT)
\end{verbatim}
}

\subsection*{Example of Yule trees with 4 leaves}
\small{
Here is a demonstration of the algorithm for the case of the Yule tree, $\alpha=\beta=0$, with $4$ leaves.

\smallskip
\begin{verbatim}
a=0; b=0; m=10000;
(bts,pts,sps)=MakePartitionAndTree(3,m,a,b)

def CountsDictWithFirstIndex(X):
    '''convert a list X into a Dictionary of counts or
    frequencies with first index of each Key saved'''
    CD = {}
    for i in range(len(X)):
        x=X[i]
        if (x in CD):
            CD[x][1] = CD[x][1]+1
        else:
            CD[x] = [i,1]
    return CD
\end{verbatim}

\smallskip
\noindent \texttt{sps} gives the list of the $m=10000$ ranked planar trees sampled by \texttt{MakePartitionAndTree}. The following function gives, for each of the trees encountered in \texttt{sps}, the theoretical probability of the tree under the Beta-splitting model with $a=b=0$ and its empirical probability (i.e., its frequency in \texttt{sps}).

\smallskip
\begin{verbatim}
BtcCnts=CountsDictWithFirstIndex(sps)
for x in BtcCnts:
    print (sps[BtcCnts[x][0]],prob_RPT(bts[BtcCnts[x][0]],a,b).N(digits=4),
            (BtcCnts[x][1]/m).N(digits=4))

((1, 3, 2), 0.1667, 0.1700)
((3, 2, 1), 0.1667, 0.1666)
((2, 1, 3), 0.1667, 0.1625)
((3, 1, 2), 0.1667, 0.1664)
((1, 2, 3), 0.1667, 0.1683)
((2, 3, 1), 0.1667, 0.1662)
\end{verbatim}

\smallskip
\noindent \texttt{bts} lists the 10000 unranked planar trees corresponding to the ranked planar trees in \texttt{sps}. The following function gives, for each of these trees, their theoretical and empirical probabilities.

\smallskip
\begin{verbatim}
BtcCnts=CountsDictWithFirstIndex(bts)
for x in BtcCnts:
    print (bts[BtcCnts[x][0]],
           prob_PT(bts[BtcCnts[x][0]],a,b).N(digits=5),(BtcCnts[x][1]/m).N(digits=5))

([., [[., .], .]], 0.16667, 0.17000)
([., [., [., .]]], 0.16667, 0.16830)
([[[., .], .], .], 0.16667, 0.16660)
([[., [., .]], .], 0.16667, 0.16620)
([[., .], [., .]], 0.33333, 0.32890)
\end{verbatim}

\smallskip
\noindent More examples can be found at \texttt{https://cloud.sagemath.com/projects/2c5f7f68-e689-4c70-a4b4-}\\
\texttt{5b5d4dc4f93f/files/2015-10-27-082849.sagews}.
}

\subsection*{A reversal property}
Although Aldous' leaf deletion property does not seem to hold in general for the random tree obtained through the generating and organizing process, at the resolution of the unranked planar trees it is possible to define a transition kernel $\overleftarrow{\mathbb{P}}$ from the set of trees with $n+1$ leaves to the set of trees with $n$ leaves in such a way that the tree obtained after $(i)$ creating a tree with $n+1$ leaves thanks to the generating and organizing process, and $(ii)$ choosing a (cherry) node to withdraw in order to come back to a tree with $n$ leaves, has the same distribution as the tree obtained from running the generating and organizing process for only $n-1$ steps. That is, writing $\mathbf{T}_n$ for the random unranked planar tree with $n$ leaves, we have for every tree $\mathbf{t}_n$ with $n$ leaves:

\begin{equation}\label{reversal}
\mathbb{P}(\mathbf{T}_n=\mathbf{t}_n)=\sum_{\mathbf{t}_{n+1}} \mathbb{P}(\mathbf{T}_{n+1}=\mathbf{t}_{n+1})\overleftarrow{\mathbb{P}}(\mathbf{t}_{n+1}\rightarrow \mathbf{t}_n).
\end{equation}

Indeed, let us set

\begin{equation}\label{proba reversal}
\overleftarrow{\mathbb{P}}(\mathbf{t}_{n+1}\rightarrow \mathbf{t}_n) = \mathbb{P}(\mathbf{T}_n=\mathbf{t}_n\, |\, \mathbf{T}_{n+1}=\mathbf{t}_{n+1}).
\end{equation}
Note that this probability is $0$ if $\mathbf{t}_n$ and $\mathbf{t}_{n+1}$ are not compatible, that is if we cannot obtain $\mathbf{t}_{n+1}$ from $\mathbf{t}_n$ by splitting one of the leaves of $\mathbf{t}_n$. Then, we trivially have

\begin{align*}
\sum_{\mathbf{t}_{n+1}} \mathbb{P}(\mathbf{T}_{n+1}=\mathbf{t}_{n+1})\overleftarrow{\mathbb{P}}(\mathbf{t}_{n+1}\rightarrow \mathbf{t}_n)& = \sum_{\mathbf{t}_{n+1}} \mathbb{P}(\mathbf{T}_{n+1}=\mathbf{t}_{n+1})\mathbb{P}(\mathbf{T}_n=\mathbf{t}_n\, |\, \mathbf{T}_{n+1}=\mathbf{t}_{n+1})\\ & = \mathbb{P}(\mathbf{T}_n=\mathbf{t}_n),
\end{align*}
which shows that \eqref{reversal} is satisfied.

Let us now give an explicit formula for the r.h.s. of \eqref{proba reversal} in the case where $\mathbf{t}_n$ and $\mathbf{t}_{n+1}$ are compatible. It is easier to come back to the resolution of ranked planar trees to compute the conditional probability appearing in the r.h.s. Indeed, as explained in the section on unranked planar trees, the probability of a given ranked planar tree $\tau_n$ does not depend on the ranking. As a consequence, conditionally on $\mathbf{T_n}=\mathbf{t}_n$, all ranked planar trees whose unranking yields $\mathbf{t}_n$ have the same probability $1/\# \mathbf{t}_n$ to be that created by the generating and organizing process, where $\# \mathbf{t}_n$ denotes the number of ranked trees corresponding to the unranked tree $\mathbf{t}_n$. Recall from the section on unranked planar trees that

$$
\# \mathbf{t}_n  = \prod_{i\in \mathcal{I}(\mathbf{t}_n)}\binom{n_i^L+n_i^R}{n_i^L}.
$$
Writing $\mathcal{T}_n$ for the random ranked planar tree with $n$ leaves and $\tau_n \prec \mathbf{t}_n$ to denote the fact that forgetting the ranking in the ranked planar tree $\tau_n$ yields $\mathbf{t}_n$, we have

\begin{align}
\mathbb{P}(\mathbf{T}_n=\mathbf{t}_n\, |\, \mathbf{T}_{n+1}=\mathbf{t}_{n+1}) & = \sum_{\tau_{n+1} \prec \mathbf{t}_{n+1}} \mathbb{P}(\mathbf{T}_n=\mathbf{t}_n\, |\, \mathcal{T}_{n+1}=\tau_{n+1}) \mathbb{P}(\mathcal{T}_{n+1}=\tau_{n+1}\, |\, \mathbf{T}_{n+1} = \mathbf{t}_{n+1}) \nonumber\\
 & = \frac{1}{\# \mathbf{t}_{n+1}} \sum_{\tau_{n+1} \prec \mathbf{t}_{n+1}} \mathbb{P}(\mathbf{T}_n=\mathbf{t}_n\, |\, \mathcal{T}_{n+1}=\tau_{n+1})\nonumber\\
 & = \frac{1}{\# \mathbf{t}_{n+1}} \sum_{\tau_{n+1} \prec \mathbf{t}_{n+1}}\sum_{\tau_n \prec \mathbf{t}_n} \mathbb{P}(\mathcal{T}_n=\tau_n\, |\, \mathcal{T}_{n+1}=\tau_{n+1}). \label{decomp}
\end{align}
Since we now work with ranked planar trees, for any tree $\tau_{n+1}$ with $n+1$ leaves the probability in the r.h.s. of \eqref{decomp} is zero unless $\tau_n$ is the tree $\tau_{n+1}^{-1}$ obtained by withdrawing the $n$-th split in $\tau_{n+1}$ (in which case the probability is $1$). Hence, the r.h.s. in \eqref{decomp} can be written

$$
\frac{1}{\# \mathbf{t}_{n+1}}\, \#\big\{\tau_{n+1}:\, \tau_{n+1}\prec \mathbf{t}_{n+1},\, \tau_{n+1}^{-1}\prec \mathbf{t}_n\big\}.
$$
But now recall that $\mathbf{t}_n$ and $\mathbf{t}_{n+1}$ are assumed to be compatible. Hence, for every tree $\tau_n$ satisfying $\tau_n \prec \mathbf{t}_n$ there is one and only one way to add a last step to obtain a tree $\tau_{n+1}\prec \mathbf{t}_{n+1}$ (namely, add the missing split in the tree and label it by $n$). As a consequence, we obtain that

\begin{equation}\label{formula reversal}
\mathbb{P}(\mathbf{T}_n=\mathbf{t}_n\, |\, \mathbf{T}_{n+1}=\mathbf{t}_{n+1}) = \frac{\# \mathbf{t}_n}{\# \mathbf{t}_{n+1}}.
\end{equation}

Note that the same definition \eqref{proba reversal} would work at the resolution of the unranked non-planar trees, but finding an explicit expression for the quantity in the r.h.s. is difficult due to the many symmetries of non-planar trees.

\begin{thebibliography}{plain}
\bibitem{mooers1997} Mooers AO, Heard SB (1997) Inferring evolutionary process from phylogenetic tree shape. \emph{Quart. Rev. Biol.} 72(1):31-54.
\bibitem{morlon2014} Morlon H (2014) Phylogenetic approaches for studying diversification. \emph{Ecology letters} 17(4):508-525.
\bibitem{colijn2014} Colijn C, Gardy J (2014) Phylogenetic tree shapes resolve disease transmission patterns. \emph{Evolution, medicine, and public health} 1:96-108.
\bibitem{colless1982} Colless DH (1982) Review of phylogenetics: The theory and practice of phylogenetic systematics. \emph{Syst. Zool.} 31:100-104.
\bibitem{sackin1972} Sackin MJ (1972) "Good" and "bad" phenograms. \emph{Syst. Zool.} 21:225-226.
\bibitem{jones2011} Jones GR (2011) Tree models for macroevolution and phylogenetic analysis. \emph{Systematic Biology} 60(6):735-746.
\bibitem{aldous2001} Aldous D (2001) Stochastic models and descriptive statistics for phylogenetic trees, from Yule to today. \emph{Statistical Science} 16(1):23-34.
\bibitem{blum2006} Blum MGB, Fran{\c{c}}ois O (2006) Which random processes describe the tree of life? A large-scale study of phylogenetic tree imbalance. \emph{Systematic Biology} 55(4):685-691.
\bibitem{yule1924} Yule GU (1924) A mathematical theory of evolution, based on the conclusions of Dr. J.C. Willis. \emph{Philos. Trans. Roy. Soc. London Ser. B} 213:21-87.
\bibitem{aldous1996} Aldous D (1996) Probability distributions on cladograms. \emph{Random discrete structures} (Springer), pp 1-18.
\bibitem{phillimore2008} Phillimore AB, Price TD (2008) Density-dependent cladogenesis in birds. \emph{PLoS Biology} 6(3):e71.
\bibitem{ford2005} Ford DJ (2005) Probabilities on cladograms: introduction to the alpha model. \emph{ArXiv preprint math/0511246}.
\bibitem{kirkpatrick1993} Kirkpatrick M, Slatkin M (1993) Searching for evolutionary patterns in the shape of a phylogenetic tree. \emph{Evolution} 47:1171-1181.
\bibitem{hagen2015} Hagen O, Hartmann K, Steel M, Stadler T (2015) Age-dependent speciation can explain the shape of empirical phylogenies. \emph{Systematic Biology} 64(3):432-440.
\bibitem{UCDMS20154} Sainudiin R, Welch D (2015) The Transmission Process. {\em UCDMS Research Report 2015/4}, School of Mathematics and Statistics, University of Canterbury, Christchurch NZ.
\bibitem{stewart2005} Stewart EJ, Madden R, Paul G, Taddei F (2005). Aging and death in an organism that reproduces by morphologically symmetric division. \emph{PLoS Biol} 3(2): e45.
\bibitem{ford2009} Ford D, Matsen FA, Stadler T (2009) A method for investigating relative timing information on phylogenetic trees. \emph{Systematic Biology} 58(2):167-183.
\bibitem{Flajolet} Flajolet P, Sedgewick R (2009) {\em Analytic Combinatorics}, Cambridge University Press, New York, NY, USA.
\bibitem{Dobrow1995} Dobrow RP, Fill JA (1995) On the Markov chain for the move-to-root rule for binary search trees. {\em The Annals of Applied Probability} 5(1):1--19.
\bibitem{CatalanCoeff2012} Sainudiin R (2012) Sequence {A}185155, {T}he {O}n-line {E}ncyclopedia of {I}nteger {S}equences, published electronically, Feb 2012.
\bibitem{Stanley1997} Stanley RP (1997) {\em Enumerative combinatorics. {V}ol. 1}, volume~49 of {\em Cambridge Studies in Advanced Mathematics}. Cambridge University Press, Cambridge, UK.
\bibitem{tajima1983} Tajima F (1983) Evolutionary relationship of DNA sequences in finite populations. \emph{Genetics} 105(2):437-460.
\bibitem{UCDMS20152} Sainudiin R, Fischer M, Cleary S, Griffiths RC (2015) Some Distributions on Finite Rooted Binary Trees. {\em UCDMS Research Report 2015/2}, School of Mathematics and Statistics, University of Canterbury, Christchurch NZ.
\end{thebibliography}
\end{document}